\documentclass[a4paper,twoside]{amsart}
\usepackage{epsf}

\title{Asymptotically efficient triangulations of the
$d$-cube}
\thanks{Research partially supported by project PB97-0358 of Spanish
{Direcci\'on General de Investigaci\'on Cient\'{\i}fica y T\'ecnica.} }

\author{David Orden}

\author{Francisco Santos}

\address{Universidad de Cantabria. Departamento de Matem\'aticas,
Estad\'{\i}stica y Computaci\'on\\
Av. Los Castros s/n, E-39005 Santander, Cantabria.
}
\email{\{ordend,santos\}@matesco.unican.es}

\date{March 2003}

\usepackage{latexsym}

\font\Bbb=msbm10 scaled 1200

\newcommand{\N}{\mbox{\Bbb N}}
\newcommand{\R}{\mbox{\Bbb R}}
\newcommand{\A}{{\mathcal A}}
\newcommand{\C}{{\mathcal C}}

\newcommand{\weight}{\operatorname{weight}}

\newcommand{\vertices}{\operatorname{vert}}
\newcommand{\conv}{\operatorname{conv}}

\newcommand{\stackldots}[1]
  {\stackrel{\hbox{\tiny $#1$}}{\ldots}}
\newcommand{\stackcdots}[1]
  {\stackrel{\hbox{\tiny $#1$}}{\cdots}}

\newtheorem{thm}{Theorem}[section]
\newtheorem{prop}[thm]{Proposition}
\newtheorem{lemma}[thm]{Lemma}
\newtheorem{coro}[thm]{Corollary}

\newtheorem{defi}[thm]{Definition}

\newtheorem{remark}[thm]{Remark}

\newenvironment{Proof} {\begin{proof}} {\end{proof}}
\begin{document}

\begin{abstract}
  Let $P$ and $Q$ be polytopes, the first of "low" dimension and the second of
"high" dimension. We show how to triangulate the product $P \times Q$
efficiently (i.e., with few simplices) starting with a given triangulation of
$Q$. Our method has a computational part, where we need to compute an
efficient
triangulation of $P \times \Delta^m$, for a (small) natural number $m$ of our
choice. $\Delta^m$ denotes the $m$-simplex.
  Our procedure can be applied to obtain (asymptotically) efficient
triangulations of the cube $I^n$: We decompose $I^n = I^k \times I^{n-k}$, for
a small $k$. Then we recursively assume we have obtained an efficient
triangulation of the second factor and use our method to triangulate the
product. The outcome is that using $k=3$ and $m=2$, we can triangulate $I^n$
with $O(0.816^{n} n!)$ simplices, instead of the $O(0.840^{n} n!)$ achievable
before.

\end{abstract}

\maketitle

\section{Introduction} 
``Simple" triangulations of the regular $d$-cube
$I^d=[0,1]^d$ have several applications, such as solving
differential
equations by finite element methods or calculating fixed points. See, for
example, \cite{Todd}. In particular, it has
brought special attention both from a theoretical point of view and from an
applied one to determine the smallest size of a triangulation of the $d$-cube
(see \cite[Section 14.5.2]{Lee} for a recent survey). Let
us point up
that the general problem of computing the smallest triangulation  of
an arbitrary
polytope is \mbox{NP-complete} even when restricted to dimension 3, see
\cite{Be-Lo-Ri}.

When we speak about
\emph{triangulations} of a polytope $P$ of dimension $d$ we mean
decompositions of
$P$ into $d$-simplices that (i) use as vertices only vertices of
$P$, and (ii) intersect face to face (i.e., forming a geometric simplicial
complex). Some authors do not require these two conditions in triangulations. 
We will always require the first one, and when the second
condition is not fulfilled, we call the decompositions 
\emph{simplicial dissections} of $P$. The \emph{size} of a
triangulation or dissection $T$ is its number of $d$-simplices
and we denote it $|T|$. It is an open question whether high dimensional
cubes admit dissections with less simplices than needed in a triangulation.
Actually, the minimum size of dissections of $I^7$ is unknown, while the
minimum triangulation is known (see below). 

The paper \cite{Lo-Ho-Sa-St} describes a general method to obtain the smallest
triangulation of a polytope $P$ as the optimal integer solution of a certain
linear program. The linear program has as many variables as
$d$-simplices with vertex set contained in the vertices of $P$ exist. That is, ${\#
\hbox{vertices}\choose {\dim(P)+1}}$ if the vertices of $P$ are 
in general position and less than that
if not. For the
$d$-cube, the direct application of this method is impossible in practice 
beyond dimension 4 or 5. With a somewhat 
similar method but simplifying the system of equations via
the symmetries of the cube, Anderson and Hughes
\cite{Anderson-Hughes} have calculated the smallest size among triangulations of
the $6$-cube and the
$7$-cube, in a computational \emph{tour-de-force} which involved  a
problem with 1\,456\,318 variables and \emph{ad hoc} ways of decomposing the
system into smaller subsystems. The smallest sizes up to
dimension 7 is shown in Table
\ref{tabla.triancubo}.

In order to compare sizes of triangulations of cubes in different dimensions,
Todd \cite{Todd} defines the \emph{efficiency} of a triangulation $T$ of the
$d$-cube to be the number $({|T|/ d!})^{1\over d}$. This number is at most 1,
since every simplex with integer vertices has a multiple of ${1/d!}$
as Euclidean
volume, so
$|T|\leq d!$. Triangulations of efficiency 1 (i.e., unimodular)
can be easily constructed in any dimension.
On the other extreme, Hadamard's inequality for determinants of matrices with
coefficients in $[-1,1]$ implies that the volume of every
$d$-simplex inscribed in the regular $d$-cube $I^d=[0,1]^d$ is at most
$(d+1)^{(d+1)/ 2}/2^d d!$. Hence, every triangulation has size at least
$2^d d!/ (d+1)^{(d+1)/ 2}$ and efficiency at least $2 / (d+1)^{d+1
\over 2d}\approx 2 / \sqrt{d+1}$.

Following the notation in \cite[Section 14.5.2]{Lee}, let 
$\phi_d$ and $\rho_d$ be the smallest size and efficiency, respectively,
of all
triangulations of the cube of dimension $d$. The number
$\phi_d$  (or some variations in which one or both of the conditions
(i) and (ii)
are not required) is known as \emph{simplexity} of the
$d$-cube. Obviously, $\rho_d=(\phi_d/d!)^{1/d}$.

In \cite{Haiman} (see also \cite{Lee}, pages 283-284), Haiman observes that a
triangulation of
$I^{k+l}$ with $t_k t_l
{k+l \choose k}$ simplices can be constructed from given
triangulations of $I^k$
and $I^l$ with $t_k$ and $t_l$ simplices respectively. With this, one easily 
concludes:

\begin{thm}[Haiman] 
\label{thm.Haiman} For every $k$ and $l$,
${\rho_{k+l}}^{k+l}\leq {\rho_{k}}^{k}{\rho_{l}}^{l}$. 
\end{thm}

\begin{coro}
\label{cor.Haiman}
The sequence
$(\rho_i)_{i\in \N}$ converges and
\[
\lim_{i\to \infty} \rho_i\le \rho_d
\qquad
\forall d\in \N.
\]
\end{coro}

\begin{proof}
Let us fix $d\in \N$ and $k\in\{1,\dots,d\}$.
Haiman's theorem implies that, for every $i\in \N$,
\[
\rho_{k+id} \le {\rho_k}^{k/(k+id)} {\rho_d}^{id/(k+id)}.
\]
Since the right-hand side converges to $\rho^d$ when $i$ grows, the $d$ subsequences of indices modulo $d$, and hence the whole sequence $(\rho_i)_{i\in \N}$, have upper limit bounded by $\rho_d$. 
An upper limit bounded by every term in the sequence must coincide with the lower limit.
\end{proof}

What is the limit of this sequence? In particular,
is it positive or is it zero? The known values of
$\rho_d$ (up to $d=7$) form a strictly decreasing sequence, as shown in Table
\ref{tabla.triancubo}, but it is not even known whether this occurs
in general.
\smallskip
\begin{table}[htb]
\begin{center}
\label{tabla.triancubo}
\begin{tabular}{| c | c | c | c | c | c | c | c | c |}
\hline Dimension & 1 & 2 & 3 & 4 & 5 & 6 & 7 & 8 \cr
\hline Smallest size ($\phi_d$)&
        1 & 2 & 5 & 16 & 67 & 308 & 1493 & $\le 11944$ \cr Smallest efficiency
($\rho_d$)&
        1 & 1 & .941 & .904 & .890 & .868 & .840 & $\le .859$ \cr Lower bound
Hadamard's ineq.&
        1 & .877 & .794 & .731 & .683 & .643 & .610 & .581  \cr Lower
bound Smith
\cite{Smith} &
        1 & 1 & .941 & .889 & .833 & .789 & .751 & .718 \cr
\hline
\end{tabular}
\caption{Smallest size and efficiency of triangulations of $I^d$.}
\end{center}
\end{table}

Concerning lower bounds, the only significant improvement to Hadamard's
inequality has been obtained in \cite{Smith}, where the same volume argument
is used, but with respect to a hyperbolic metric. The last row of the Table
\ref{tabla.triancubo} shows the lower bound obtained, translated into
efficiency of triangulations. For small dimensions, it is an excellent
approximation of the smallest efficiency. Asymptotically, it only increases
the lower bound obtained using Hadamard's inequality by a constant factor
$\sqrt{3/ 2}$.

In this paper we propose a method to obtain efficient
triangulations of a product polytope $P\times Q$ starting from a
triangulation of $Q$ and another of $P\times \Delta^{m-1}$, where
$\Delta^{m-1}$ denotes a simplex of dimension $m-1$ and $m$ is any
relatively small number. We apply this with $P$ being a small-dimensional cube
and $Q$ a high-dimensional one, iteratively. This allows to obtain
asymptotically efficient triangulations of arbitrarily high-dimensional cubes
from any (efficient) triangulation of $I^l\times \Delta^{m-1}$. 
Sections 3, 4, and 5 explain our method, which is first outlined in
Section 2. If the reader is happy with dissections,
we cannot offer better efficiencies for them than for triangulations, but
at least he or she can skip Section 5, which contains most of the
technicalities in this paper.

The asymptotic efficiency of the triangulations obtained, clearly,
depends on how good our triangulation of  $I^l\times \Delta^{m-1}$
is. Finding the triangulation of $I^l\times\Delta^{m-1}$ which is optimal for
our purposes reduces to an integer programming problem, similar to 
finding the smallest triangulation of that polytope (actually, it is the same
system of linear
equations, with a different objective function). Using the linear
programming software {\tt CPLEX}, we have solved the system for some values of
$l$ and $m$. The best triangulation we have found is one of 
$I^3\times\Delta^{2}$, with which we obtain
\[
\lim_{i\to \infty} \rho_i\le \sqrt[3]{{44/3}\over 27}\approx 0.8159
\]

The best bound existing before was $\lim_{i\to \infty} \rho_i\le
\rho_7 =0.840$. In other words, 
we prove that (asymptotically) the $d$-cube can be
triangulated with $0.8159^d d!$ simplices, instead of the
$0.840^d d!$ achievable before.

It has to be observed that, even if the particular triangulation of $I^3\times\Delta^{2}$ that we use
was obtained by an intensive computer calculation, once the triangulation is found it is a simple task 
to check that it is indeed a triangulation and compute the asymptotic efficiency obtained from it.
Actually, in Section \ref{sec.rholm} we use the so-called Cayley Trick \cite{Hu-Ra-Sa}
to do this checking with no
need of computers at all.  We also use the Cayley Trick
to explore the minimum efficiency that can be obtained from the product $I^2\times\Delta^{k}$ for any $k$.
We briefly explain the trick in Section
\ref{sec.cayley}, where we also interpret our whole construction in terms of it.

\section{Overview of the method and results}
\label{sec.overview}

We start describing Haiman's proof of Theorem \ref{thm.Haiman}, which is related to our method.
Let $T_k$ and $T_l$ be triangulations of the regular cubes
$I^k$ and $I^l$, respectively. The product $T_k\times T_l$
of the two triangulations gives 
a decomposition of the cube $I^k\times I^l =I^{k+l}$ into $|T_k|\cdot|T_l|$
subpolytopes, each of them isomorphic to the product of simplices
$\Delta^k\times\Delta^l$. It is well
known that every triangulation of
$\Delta^k\times\Delta^l$ has size ${k+l}\choose{k}$
(see, e.g., \cite[Chapter 7]{GKZbook}). Hence, 
refining $T_k\times T_l$ in an arbitrary way one gets a triangulation of
$I^{k+l}$ of size $|T_k|\cdot|T_l|\cdot {{k+l}\choose{k}}$. The
implication of this is that starting with a triangulation of a cube
$I^k$ of a certain efficiency $\rho$, one can construct a sequence of 
triangulations of $I^{nk}$ for $n\in \N$ whose asymptotic
efficiencies converge to $\rho$.

Our method is, in a way, similar.
Starting from a triangulation of
$I^{n-1}$ and another one of $I^l\times\Delta^{m-1}$, we get one of 
$I^{l+n-1}$ with the
following general method to triangulate $P\times Q$ starting from a
triangulation of $Q$ and another one of $P\times\Delta^{m-1}$:

\begin{enumerate}
\item[(1)] We first show (Sections \ref{sec.subdivision} and
\ref{sec.refinamiento}) how to obtain triangulations of
$P\times\Delta^{n-1}$ from triangulations of
$P\times\Delta^{m-1}$, where $n-1=\dim (Q)$ is supposed to be much bigger
than $m-1$.
We call our triangulations \emph{multi-staircase triangulations}.

\item[(2)] A triangulation of $Q$ induces, as in Haiman's method, a
decomposition of $P\times Q$ into polytopes isomorphic to
$P\times\Delta^{n-1}$. Each of them can be triangulated using
the previous paragraph, although this in principle only 
gives a dissection of $P\times Q$; if we want a
triangulation, we have to apply (1) to all the polytopes
$P\times\Delta^{n-1}$ in a  compatible way.
In Section \ref{sec.triang}
we will show how to do this using an $m$-coloring of the vertices of
$Q$. 

\item[(3)]
The analysis of the size of the triangulation obtained will also be carried
out in Section \ref{sec.triang}. 
\end{enumerate}

In step (2), the final size of the
triangulation is just the sum of the individual sizes of the triangulations
used for the different subpolytopes $P\times \Delta^{n-1}$. In particular, if
we are just interested in obtaining dissections, we can take an efficient
triangulation of $P\times \Delta^{n-1}$ and repeat it in every
subpolytope.  

In step (1) the computation of the size is more complicated and to state it in
a simple way we introduce the following definitions. When we speak about
simplices in a polytope $P$ we will implicitly suppose that its vertices are
vertices of $P$, and we will identify the simplex with its vertex set.

\begin{defi}
\label{def.main}
\rm Let $P$ be a polytope of dimension $l$ and let 
$\tau$ be a simplex of dimension $l+m-1$ in $P\times
\Delta^{m-1}$. Let
$\{v_1,\ldots, v_m\}$ be the vertices of $\Delta^{m-1}$. Then,
$\tau$, understood as a vertex set, decomposes as
$\tau=\tau_1\cup \cdots \cup\tau_m$ with $\tau_i\subset P\times\{v_i\}$.
\begin{enumerate}
\item[(i)] We define the sequence $(|\tau_1|-1,\ldots,|\tau_m|-1)$ to be the
\emph{type} of the simplex $\tau$.   The \emph{weight} of a simplex
$\tau$ of type
$(t_1,\ldots,t_m)$ is the number
$ { 1 \over \prod_{i=1}^m t_i !}.
$
\item[(ii)] The \emph{weighted size} of a triangulation $T$ of
$P\times
\Delta^{m-1}$ is $\sum_{\tau\in T} \weight(\tau)$, and the
\emph{weighted efficiency} is
$$
\sqrt[l]{ \sum_{\tau\in T} \weight(\tau)\over m^l}
$$
\end{enumerate}
\end{defi}

With this, our main result can be stated as:

\begin{thm}
\label{thm.expectedsize}
Consider polytopes $P$ of dimension $l$ and $Q$ of dimension $n-1$. Let $m$
be such that 
$m\leq n$. Given a triangulation $T_0$ of $P\times\Delta^{m-1}$ of
weighted size $t_0$  and a triangulation $T_Q$ of $Q$, then there are
triangulations  of $P\times Q$ with size at most
\[
|T_Q|t_0\left({n\over{m}}+l\right)^l.
\]
\end{thm}

Using this method to triangulate $I^{l+n-1}=I^l\times I^{n-1}$ starting from a
triangulation of $I^{n-1}$ and another one of $I^l\times\Delta^{m-1}$ we will
conclude:

\begin{thm}
\label{thm.triangulaciones} If there is a triangulation $T_0$ of
$I^l\times\Delta^{m-1}$ with weighted efficiency $\rho_0$, then
for every $\epsilon>0$ and all $n$ bigger than 
${lm\rho_0/\epsilon}$ we have
\[
   \quad {\rho_{n+l-1}}^{n+l-1}\leq{\rho_{n-1}}^{n-1}{(\rho_{0}+\epsilon)}^{l}.
\]
As a consequence,
\[
\lim_{i\to \infty} \rho_i\le \rho_0.
\]
\end{thm}

Observe that the definition of weighted efficiency makes
sense even in the case $m=1$, where all the simplices of
$P\times\Delta^{0}$ have the same type, equal to $(l)$ if $\dim(P)=l$, and the
same weight, equal to $1/l!$. In particular, the weighted 
efficiency of a triangulation of $I^l\times \Delta^{0}$ is the same as the
usual ``non-weighted" one. Another common point between efficiency and
weighted efficiency is that the weighted
efficiency of a triangulation of
$I^l\times \Delta^{m-1}$ is always less or equal to 1, and it is 1 if and only
if the triangulation is unimodular; i.e., if every simplex has volume
$1/(l+m-1)!$ (this is proved in Section \ref{sec.cayley}).
\label{futurdem.efftriang1}

We now describe the practical results obtained. The last theorem leads us to
study the smallest weighted efficiency of triangulations of $I^l\times
\Delta^{m-1}$; let us denote it by $\rho_{l,m}$. This number can be calculated
minimizing a linear form over the so-called \emph{universal polytope} of all
the triangulations of $I^l\times \Delta^{m-1}$. 

The definition of this
universal polytope for triangulations of an arbitrary polytope $P$ is as
follows (see \cite{Lo-Ho-Sa-St} for further details): 
let $\Sigma(P)$ be the set of all the simplices of maximal dimension which 
use as vertices only vertices of $P$. Given a triangulation $T$ of
$P$, its \emph{incidence vector} $V_T\in\R^{\Sigma(P)}$ has a 1 in the
coordinates corresponding to simplices of $T$ and a 0 in the others. The
\emph{universal polytope} of $P$ is
$\conv\{V_T : T \mbox{ is a triangulation of }P\}$.

In our case, to calculate the smallest weighted size (and efficiency), we have to minimize
over the universal polytope of \mbox{$I^l\times \Delta^{m-1}$} 
the linear form having as coefficient of each
simplex its weight ${1/{\prod_{i=1}^{m}t_i!}}$. (To minimize non-weighted efficiency and size
one just uses the functional with all
coefficients equal to 1). One of the results in
\cite{Lo-Ho-Sa-St} is a description of the vertices of the universal
polytope as
integer solutions of a certain system of linear inequalities derived from the
oriented matroid of $P$. Thus, our minimization problem is restated
as an integer
linear programming problem.

In order to apply our method, we have used the program {\tt UNIVERSAL
BUILDER} by
Jes\'us A. de Loera and Samuel Peterson \cite{Universal}. 
Given as input the vertices
of $P$, the
program generates the linear system of equations defining the
universal polytope
of $P$. The output is a file readable by the linear programming software {\tt
CPLEX}. We have created a routine that generates our particular objective
function. Table
\ref{tabla.eficiencias} shows the results obtained in the cases we
have been able
to solve. Note that $\rho_{l,1}=\rho_l$ and so the column $m=1$ in Table
\ref{tabla.eficiencias}  is the same as
the second row of Table \ref{tabla.triancubo}.

The fact that $\rho_{1,m}=1$ for every $m$ reflects that every triangulation of
the prism $I\times
\Delta^{m-1}$ is unimodular. In the case of $\rho_{2,m}$ we prove (Subsection
\ref{subsec.l=2}) that the smallest weighted efficiency is always
$\sqrt{\lceil 3m^2/4 \rceil/m^2}$.  That is to say, 
$\sqrt{3/4}$ for even $m$ and $\sqrt{3/4}+\Theta(m^{-2})$ for odd $m$.

\label{rhos.cplex}
The computation of $\rho_{3,3}$ involved a
system with 74\,400 variables, whose resolution by {\tt{CPLEX}} took
37 hours of CPU on a SUN UltraSparc.
\smallskip
\begin{table}[htb]
\begin{center}
\begin{tabular}{| c | c | c | c | c |}
\hline
  $l\setminus m$ & 1 & 2 & 3 & $\ge 3$ \\
\hline 1 & 1 & 1 & 1 & 1 \cr
\hline 2 & 1 &
$\sqrt{3\over 4} \approx 0.866$ & $\sqrt{7\over 9}\approx 0.8819$ &
$\ge\sqrt{3\over 4}$ \\
\hline 3 & $\sqrt[3]{{5/6}\over{1}}\approx 0.941$ &
$\sqrt[3]{{14/3}\over 8}\approx 0.8355$ &
$\sqrt[3]{{44/3}\over 27}\approx 0.8159$ & \\
\hline
\end{tabular}
\end{center}
\caption{Values of $\rho_{l,m}$ for $l\le 2$ or for $l=3$ and $m\le 3$.}
\label{tabla.eficiencias}
\end{table}

\section{Polyhedral subdivision of $P\times\Delta^{k_1+\cdots
+k_m-1}$ induced by
a  polyhedral subdivision of $P\times\Delta^{m-1}$}
\label{sec.subdivision}

We call \emph{polyhedral subdivisions} of a polytope $P$ its face-to-face 
partitions into subpolytopes which only use vertices of $P$ as vertices.

Let $P$ be a polytope of dimension $l$. Let $m$ and
$k_1,\ldots, k_m$ be natural numbers and let us call $n:=k_1+\cdots +k_m$. Let
$v_1,\ldots,v_m$ be the $m$ vertices of the standard simplex
$\Delta^{m-1}$ and
$v_1^1,\ldots,v_{k_1}^1,\ldots\ldots,v_1^m,\ldots,v_{k_m}^m$ the vertices of
$\Delta^{n-1}$. Observe that, implicitly, we have the following 
surjective map:

\begin{center}
\begin{tabular}{ccc}
$\vertices(\Delta^{n-1})$ &$\to$ &$\vertices(\Delta^{m-1})$  \cr
                  $v_j^i$ &$\mapsto$ & $ v_i$  \cr
\end{tabular}
\end{center} 
This map uniquely extends to an affine projection
$\pi_0:\Delta^{n-1} \to \Delta^{m-1}$. In turn, this induces a projection
\begin{center}
\begin{tabular}{ccc}
$\pi={\Bbb 1}\times \pi_0:
       P\times\Delta^{n-1}$ &$\to$ &$ P \times\Delta^{m-1}$  \cr
        \hspace{1.5cm}$(p,a)$ &$\mapsto$ & $ (p,\pi_0(a))$  \cr
\end{tabular}
\end{center}

Given the projection $\pi$ and a polyhedral subdivision $S$ of the target
polytope $P \times\Delta^{m-1}$, it is obvious that the inverse images
$\pi^{-1}(B)$ of the subpolytopes of $S$ form a subdivision of
$P\times\Delta^{n-1}$ into subpolytopes matching face to face.
In a more general projection it would not be true 
that those subpolytopes use as vertices only vertices of
$P\times\Delta^{n-1}$. But it is true in our case:

\begin{lemma}
\label{lema.proyeccion} Let $B\subset P\times \Delta^{m-1}$ be a
subpolytope with
$\vertices(B)\subset\vertices(P\times\Delta^{m-1})$. Let $\pi$ be the
projection
considered before. Let be $\tilde B=\{(p,v^i_j) : (p,v_i)\in \vertices(B)\}$.
Then:
\[
\pi^{-1}(B)\subset\conv(\tilde B).
\]
\end{lemma}

\begin{Proof} Let $(p,a)$ be a point of $\pi^{-1}(B)$, so
$(p,\pi_0(a)) \in B$. Let us write $a$ as a convex combination of the
vertices of
$\Delta^{n-1}$, that is $a=\sum_{i=1}^m \sum_{j=1}^{k_i}
\lambda_j^i  v_j^i$, with
$\lambda_j^i\ge 0$, $\forall i,j$ and $\sum_i\sum_j \lambda_j^i=1$.
Without loss
of generality, we can suppose that none of the sums
$\sum_{j=1}^{k_i} \lambda_j^i$ is zero, otherwise everything
``happens" on a face
of $\Delta^{m-1}$ and we can restrict the statement to that face.

Let $(p_1^1,v_1),\ldots,(p_{l_1}^1,v_1)$,
$(p_1^2,v_2),\ldots,(p_{l_2}^2,v_2)$,
$\ldots$,
$(p_1^m,v_m),\ldots,(p_{l_m}^m,v_m)$ be the vertices of $B$, so we have
$
\tilde{B}= \{ (p_h^i,v_j^i) : i=1,\ldots,m, j=1,\ldots,k_i, h=1,\ldots,l_i\}.
$ We write now $(p,\pi_0(a))$ as convex combination of the vertices of
$B$, that is,
\[ (p,\pi_0(a))= \sum_{i=1}^m \sum_{h=1}^{l_i} \mu_h^i \,(p_h ^i,v_i),
\] with $\mu_h^i\ge 0$, $\forall i,h$, and $\sum_{i} \sum_{h}
\mu_h^i=1$. Observe that $\sum_{j=1}^{k_i}\lambda_j^i
=\sum_{h=1}^{l_i} \mu_h^i$
for every $i$, because
$\pi_0(a)=\sum_{i=1}^m \sum_{j=1}^{k_i}\lambda_j^i  v_i$. Then, it is
easy to check
that:
\[ (p,a)=
\sum_{i=1}^m \sum_{j=1}^{k_i} \sum_{h=1}^{l_i} {\lambda_j^i
\mu_h^i\over \sum_{j=1}^{k_i} \lambda_j^i}  (p_h ^i,v_j^i)
\mbox{\hskip1cm and \hskip1cm} 1=\sum_{i=1}^m \sum_{j=1}^{k_i}
\sum_{h=1}^{l_i}  {\lambda_j^i \mu_h^i\over \sum_{j=1}^{k_i}
\lambda_j^i}
\] That is, $(p,a)$ is a convex combination of points of $\tilde B$,
as we wanted
to prove.
\end{Proof}

\begin{coro} Under the previous conditions;
\begin{enumerate}
\item[(i)] $\pi^{-1}(B)=\conv(\pi^{-1}(\vertices(B))).$
\item[(ii)] $\vertices(\pi^{-1}(B))=\tilde{B}$.
\end{enumerate}
\end{coro}

\begin{Proof} In both equalities, the inclusion from right to left is
trivial. In
the second one, observe that if
$(p,v_j^i)$ is in $\tilde{B}$, then $(p,v_i)$ is a vertex of $B$. Therefore,
$(p,v_j^i)\in\pi^{-1}(B)$ and it is a vertex of $P\times
\Delta^{n-1}$, which implies that it is a vertex of $\pi^{-1}(B)$.

Inclusions from left to right follow from  Lemma \ref{lema.proyeccion}, because:

\noindent (1) $\tilde{B}\subset \pi^{-1}(\vertices(B)) \Rightarrow
\pi^{-1}(B)\subset \conv(\tilde{B})\subset
\conv(\pi^{-1}(\vertices(B))).$

\noindent (2) $\pi^{-1}(B)\subset \conv(\tilde{B}) \Rightarrow
\vertices(\pi^{-1}(B)) \subset
\vertices(\conv(\tilde{B}))=\tilde{B}$, where the
last equality follows from the fact that all the elements of
$\tilde{B}$ are vertices of $P\times
\Delta^{n-1}$.
\end{Proof}

\begin{coro} 
\label{coro.subdivision}
Every polyhedral subdivision $S$ of $P\times
\Delta^{m-1}$ induces a polyhedral subdivision
$\tilde S:=\pi^{-1}(S)=\{\pi^{-1}(B) : B\in S\}$ of $P\times \Delta^{n-1}$.
Furthermore, the vertices of each subpolytope $\pi^{-1}(B)$ in this
subdivision are  $\tilde B:=\{(p, v^i_j) : (p,v_i)\in\vertices (B) \}$.
\qed
\end{coro}

\section{Triangulation of $P\times\Delta^{k_1+\cdots +k_m-1}$ induced by a
triangulation of $P\times\Delta^{m-1}$}
\label{sec.refinamiento}

We will suppose now that the polyhedral subdivision $S$ of $P\times
\Delta^{m-1}$ is a triangulation. A convenient graphic representation of the
vertices of the polytope $P\times \Delta^{m-1}$ is as a grid whose rows
represent vertices of $P$ and whose $m$ columns represent the vertices
$v_1,\ldots,v_m$ of $\Delta^{m-1}$. In order to represent a subset of vertices
of $P\times \Delta^{m-1}$ we just mark the corresponding squares. In the same
way we can represent the vertices of $P\times \Delta^{k_1+\cdots +k_m-1}$, but
now it is convenient to divide the grid horizontally in blocks, each of them
corresponding to each $k_i$ and containing the vertices $v_1^i,\ldots
,v_{k_i}^i$ of $\Delta^{k_1+\cdots +k_m-1}$.

\begin{figure}[htb]
\begin{center}
         \epsfxsize=5 in\leavevmode \epsfbox{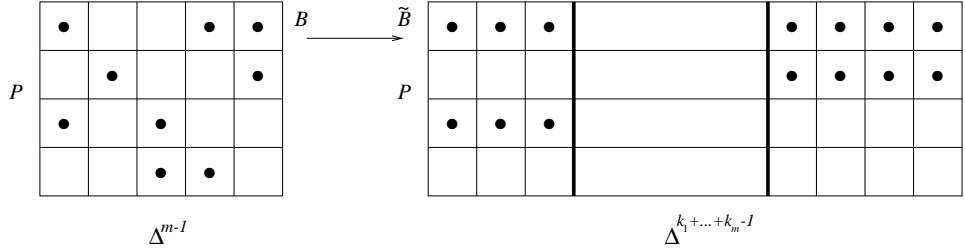}
\end{center}
\caption{How to obtain $\tilde B$ from $B$.}
\label{fig.bloques}
\end{figure}

Figure \ref{fig.bloques} shows how to obtain, with the notation of the previous
section, the set $\tilde B$ associated to a simplex $B\in S$ in this graphic
representation. In the $i$-th block, rows corresponding to vertices
$(p,v_i)$ in
$B$ have all its squares marked. Restricting $\tilde B$ to that block gives
precisely
\mbox{$\conv(\{p_1^i,\ldots,p_{l_i}^i\})\times
\conv(\{v_1^i,\ldots,v_{k_i}^i\})$}, with the notation used for the vertices of
$B$ in the proof of Lemma
\ref{lema.proyeccion}.

\begin{figure}[htb]
\begin{center}
         \epsfxsize=3 in\leavevmode \epsfbox{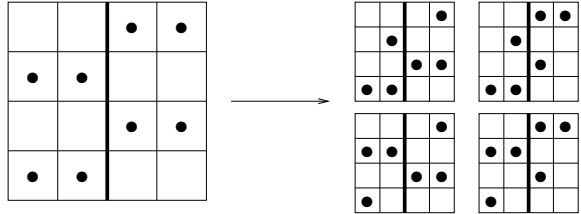}
\end{center}
\caption{The four multi-staircases (right) forming the multi-staircase
triangulation of the polytope $\tilde B$ in the left.}
\label{fig.multiescaleras}
\end{figure}

Since $B$ is a simplex, any subset of its vertices forms also a simplex; thus,
the restriction of $\tilde B$ to each block is a product of two
simplices. Let us
recall that the \emph{staircase triangulation} of the product of two simplices
$\Delta^k\times\Delta^l$ is the one whose $k+l\choose k$ simplices are all the
possible monotone  staircases in a grid of size $(k+1)\times (l+1)$. Following
this analogy, we define \emph{multi-staircases} as follows:

\begin{defi}
\rm Let $\tilde B = \pi^{-1}(B)$ be a subpolytope of 
$P\times \Delta^{k_1+\cdots +k_m-1}$, of the kind obtained in Section
\ref{sec.subdivision}. 
\begin{enumerate}
\item[(i)] A \emph{multi-staircase} in $\tilde B$ is any subset of
vertices which restricted to every block forms a monotone staircase.
\item[(ii)] The \emph{multi-staircase triangulation} of $\tilde B$ is the one
which has as simplices the different multi-staircases (see Figure
\ref{fig.multiescaleras}).
\end{enumerate}
\end{defi}

\begin{lemma} The multi-staircases indeed
form a triangulation of $\tilde B$ and taking
the multi-staircase triangulations of the different
$\tilde B$'s obtained from a triangulation of $P\times \Delta^{m-1}$
we get a triangulation of $P\times \Delta^{k_1+\cdots + k_m -1}$, which we
call \emph{multi-staircase triangulation}.
\end{lemma}

\begin{figure}[htb]
\begin{center}
         \epsfxsize=3 in\leavevmode \epsfbox{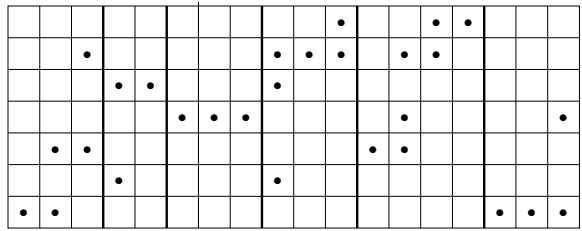}
\end{center}
\caption{An example of a multi-staircase.}
\label{fig.multiescal}
\end{figure}

\begin{Proof} 
It is clear that multi-staircases form full-dimensional simplices in
$\tilde{B}$. A way to prove that a collection $T$ of full-dimensional
simplices in a polytope $\tilde{B}$ is a triangulation is to show that:
\begin{enumerate}
\item They induce a triangulation on one face of $\tilde{B}$.

\item For each full-dimensional simplex in the collection, the removal
of any single vertex produces a codimension one simplex either
\begin{itemize}
\item lying on a facet of $\tilde{B}$ and not contained in any other simplex
of $T$, or
\item contained in exactly another full-dimensional
simplex of $T$ which is separated from the first one by their common facet.
\end{itemize}
\end{enumerate}

In our case, the first condition follows by induction
on $k_1+\cdots+k_m$, the base case being $k_1+\cdots+k_m=m$. 

For the second condition, let $\sigma$ be a multi-staircase, and let
$(p,v)$ be a vertex in it. Then, one of the following three things
happens (Figure \ref{fig.multiescal} gives an example of a multi-staircase):
\begin{itemize}
\item If $(p,v)$ is the only point of $\sigma$ in its column, 
then no other multi-staircase contains $\sigma
\setminus\{(p,v)\}$. Removing $(p,v)$ produces indeed a
codimension one simplex contained in a facet
$P\times\Delta^{k_1+\cdots+k_m-2}$ of
$P\times\Delta^{k_1+\cdots+k_m-1}$.

\item If $(p,v)$ is the only point of $\sigma$ in a row within a block,
then no other multi-staircase contains $\sigma
\setminus\{(p,v)\}$. Removing $(p,v)$ produces indeed a
codimension one simplex in a facet of $\tilde{B}$ of the form
$\pi^{-1} (B\setminus \{(p,\pi_0(v))\})$ (remember that
$B=\pi(\tilde{B})$ is a simplex).

\item If $(p,v)$ is an elbow in the multi-staircase, then removing it
leads to a unique different way of completing the
multi-staircase. More precisely, let $(p',v)$ and $(p,v')$ be the
points of the multi-staircase adjacent to $(p,v)$. Then removing
$(p,v)$ and inserting $(p',v')$ produces the other possible
multi-staircase. The obvious affine dependency
$(p,v) + (p',v')=(p,v') + (p',v)$ implies that the two multi-staircases
lie in opposite sides of their common facet.
\end{itemize}

\end{Proof}

\begin{lemma}
\label{lema.size} Let $l_i$ be the number of vertices of $B$ in the
$i$-th column. Then the multi-staircase triangulation of $\tilde B$ has exactly
$\prod_{i=1}^{m} {k_i-1+l_i-1 \choose k_i-1}$ simplices.
\end{lemma}

\begin{Proof} In each sub-block there are ${k_i-1+l_i-1 \choose
k_i-1}$ possible monotone staircases.
\end{Proof}

\section{A triangulation of $P\times Q$}
\label{sec.triang}

We will consider polytopes $P$ of dimension $l$ and $Q$ of dimension $n-1$. We
assume we are given a triangulation
$T_Q$ of
$Q$, which induces a decomposition of
$P\times Q$ into cells isomorphic to $P\times \Delta^{n-1}$. Then, any
decomposition
$n:=k_1+\cdots +k_m$ allows us to apply the procedure of Sections
\ref{sec.subdivision} and  \ref{sec.refinamiento} to triangulate the
cells $P\times \Delta^{n-1}$ starting from a triangulation $T_0$ of 
$P\times\Delta^{m-1}$.

There are two important tasks remaining: First, show that the triangulations of
the different $P\times
\Delta^{n-1}$ can be achieved in a coordinated way to obtain a real
triangulation
of $P\times Q$; second, analyze the efficiency of that triangulation. Both of
them will be done in this section, using the following trick:

Consider a partition of the vertices of $Q$ into $m$ ``colors". Then, in
each subpolytope $P\times
\Delta^{n-1}$ of $P\times Q$ the vertices of the factor $\Delta^{n-1}$ are
colored as well. We use this coloring to construct the projections
$\Delta^{n-1}\to\Delta^{m-1}$ we need for each of them. Then, on each
common face
(isomorphic to $P\times\Delta^{k}, k<n-1$) of two of the cells $P\times
\Delta^{n-1}$ we get the multi-staircase triangulation induced by the
$m$-coloring of
$Q$ restricted to that face. Hence, the triangulations of the different
cells $P\times \Delta^{n-1}$ intersect face-to-face, 
and we obtain a triangulation
$T_{P\times Q}$ of $P\times Q$, which we call \emph{multi-staircase
triangulation} of $P\times Q$.

In order to analyze the size of $T_{P\times Q}$, we will
suppose that
the $m$-coloring of the vertices of $Q$ is chosen at random with a
uniform distribution.

For each $\sigma\in T_Q$ let $\sigma_i:=\{\mbox{vertices of
$\sigma$ colored $i$}\}$. And for each
$\tau\in T_0$, $\tau_i:=\{\mbox{vertices of $\tau$ over the $i$-th vertex of
$\Delta^{m-1}$}\}$. By Lemma \ref{lema.size}, the triangulation $T_{P\times Q}$
we obtain has
$$
\sum_{\sigma\in T_Q}
\sum_{\tau\in T_0}
\prod_{i=1}^{m} {|\sigma_i|-1+|\tau_i|-1 \choose |\tau_i|-1}
$$ simplices. The expected value of this sum, when the coloring is random,
equals the sum of the expected values. So let us fix a pair of
simplices $\tau\in
T_0$ and $\sigma\in T_Q$ and calculate the expected value of
\[
\prod_{i=1}^{m} {|\sigma_i|-1+|\tau_i|-1 \choose |\tau_i|-1}
\] 
We call $l_i:=|\tau_i|$ and $k_i:=|\sigma_i|$. The $l_i$'s are considered
constants, while the $k_i$'s are random variables depending on the coloring.
They follow a multinomial distribution, in which the probability of the
$m$-tuple $(k_1,\ldots,k_m)$ is $ P(k_1,\ldots,k_m)= {{n!}\over{k_1!\cdots
    k_m!}{m^n}}.  $

Since ${{x}\choose{n}}={x^{\underline{n}}/{n!}}$, where
$x^{\underline{n}}:=x(x-1)\cdots(x-n+1)$ is the $n$-th
\emph{falling power} of $x$, we can write:
\[
E(\prod_{i=1}^{m}{k_i-1+l_i-1 \choose l_i-1})=
\frac{E(\prod_{i=1}^{m}(k_i-1+l_i-1)^{\underline{l_i-1}})}
     {\prod_{i=1}^{m} (l_i-1)!}
\]
For the numerator in the RHS of this equation we will use the following result
from \cite{Whittle} (Theorem 4.4.4):
%\note{more precise reference}

\begin{thm}
\label{thm.libroWhittle} Let $X_1,\ldots,X_m$ be scalar random variables. Then,
the expected value of
$\prod_{i=1}^{m}X_i^{\underline{v_i}}$ is given by the formal calculus
\[ E(\prod_{i=1}^{m}X_i^{\underline{v_i}})=
\left[
\begin{array}{c} {{\partial^{v_1+\cdots+v_m}}
\over{\partial z_1^{v_1}\cdots\partial z_m^{v_m}}} E(\prod_{i=1}^{m}z_i^{X_i})
\end{array}
\right]_{z_1=\cdots =z_m=1}
\] in which extra variables $z_i,\ i=1,\ldots,m$ appear with formal purposes.
\end{thm}

\begin{lemma}
\label{lema.expectation} Let $l_1,\ldots,l_m$ be positive integers 
with $\sum l_i=l+m$ and
let $k_1,\ldots,k_m$ be random variables  obeying a multinomial distribution
with $\sum k_i = n$. Then,
\[ E(\prod_{i=1}^{m}(k_i-1+l_i-1)^{\underline{l_i-1}})\leq
\left(l+{n\over m}\right)^{l}
\]
\end{lemma}

\begin{Proof} If some $l_i$ equals 1 we can neglect it in the statement: we
  remove it and will still have $\sum (l_i-1)=l$. Hence we will assume $l_i\ge
  2$ for every $i$.

We will use the following formula for expectations under a multinomial
distribution, again taken from \cite{Whittle} (p. 53, last formula in the proof of Theorem 4.2.1):
\[ 
E(\prod_{i=1}^{m}z_i^{k_i})=
\left(
\sum_{i=1}^{m}{1\over m}z_i
\right)^{n}
\]

Since
$\sum_{i=1}^{m}(l_i-1)=l$ and only the $k_i$ depend on the random process,
applying Theorem \ref{thm.libroWhittle} to the random variables
$X_i:=k_i+l_i-2$ and using the equality above we get: 
\[ 
E(\prod_{i=1}^{m}(k_i-1+l_i-1)^{\underline{l_i-1}})=
\left[
\begin{array}{c} {{\partial^{l}}
\over{\partial z_1^{l_1-1}\cdots\partial z_m^{l_m-1}}}
E(\prod_{i=1}^{m}z_i^{k_i-1+l_i-1})
\end{array}
\right]_{\overline{z}=1}=
\]
\[ =
\left[
\begin{array}{c} {{\partial^{l}}
\over{\partial z_1^{l_1-1}\cdots\partial z_m^{l_m-1}}}
\left(
\prod_{i=1}^{m}z_i^{l_i-2}E(\prod_{i=1}^{m}z_i^{k_i})
\right)
\end{array}
\right]_{\overline{z}=1}=
\]
\[ =
\left[
\begin{array}{c} {{\partial^{l}}
\over{\partial z_1^{l_1-1}\cdots\partial z_m^{l_m-1}}}
\left(
\prod_{i=1}^{m}z_i^{l_i-2}\left(
\sum_{i=1}^{m}{1\over m}z_i
\right)^{n}
\right)
\end{array}
\right]_{\overline{z}=1}\stackrel{(\ast)}{\leq}
\prod_{i=1}^{m}\left(l_i-2+{n\over m}\right)^{l_i-1}
\leq
\]
\[
\leq\prod_{i=1}^{m}\left(l+{n\over m}\right)^{l_i-1}=
\left(l+{n\over m}\right)^{l},
\] 
where the last inequality comes from $l_i\leq l+1$ (since $\sum (l_i-1)=l$
and $l_i-1\geq 1$), and the one marked with an asterisk needs to be proved.
For this we call
$F_{l_1,\ldots,l_m}(\overline{z}):= \prod_{i=1}^{m}z_{i}^{l_i-2}$
and $G_n(\overline{z}):=\left(\sum_{i=1}^{m}{1\over m}z_i\right)^n$.
Using that 
\[
{\partial^k \over \partial x^k} f(x) g(x)= \sum_{j=0}^k {k\choose j}
{\partial^j \over \partial x^j} f(x) {\partial^{k-j} \over \partial x^{k-j}}
g(x)
\]
we come up to
\[ {{\partial^{l_i-1}}\over{\partial z_i^{l_i-1}}} 
\left(
F_{l_1,\ldots,l_m}(\overline{z}) G_{n}(\overline{z})
\right)
=
\sum_{j=0}^{l_i-1}{l_i-1\choose j}(l_i-2)^{\underline{l_i-1-j}}\left({1\over
m}\right)^{j}\!\!\!n^{\underline{j}}
\ F_{l_1,\ldots,j+1,\ldots,l_m}(\overline{z})G_{n-{j}}(\overline{z}).
\]

\noindent Therefore, since
\[
{{\partial^{l}}\over{\partial {z_1}^{l_1-1}\cdots\partial {z_m}^{l_m-1}}}
     h(\overline{z})= 
{{\partial^{l_m-1}}\over{\partial
{z_m}^{l_m-1}}}\left(\cdots\left({{\partial^{l_1-1}}\over{\partial
{z_1}^{l_1-1}}}
\left(h(\overline{z})\right)\right)\cdots\right)
\]
and
\[
\left[F_{h_1,\ldots,h_m}(\overline{z})\right]_{\overline{z}=1}=
\left[G_{t}(\overline{z})\right]_{\overline{z}=1}=1,
\]
we get
\[
\left[ {{\partial^{l}}\over{\partial z_1^{l_1-1}\cdots\partial z_m^{l_m-1}}}
F_{l_1,\ldots,l_m}(\overline{z})G_{n}(\overline{z})
\right]_{\overline{z}=1} =
\]
\[
=
\sum_{j_1=0}^{l_1-1}  \cdots  \sum_{j_m=0}^{l_m-1} 
\left(
\prod_{i=1}^m
{l_i-1\choose j_i}(l_i-2)^{\underline{l_i-1-j_i}}
\left({1\over m} \right)^{j_i}
(n-j_1-\cdots -j_{i-1})^{\underline{j_i}}
\right)
=
\]
\[
=
\sum_{j_1=0}^{l_1-1}  \cdots  \sum_{j_m=0}^{l_m-1} 
\left(
\prod_{i=1}^m
{l_i-1\choose j_i}(l_i-2)^{\underline{l_i-1-j_i}}
\left({1\over m} \right)^{j_i}
\right)
n^{\underline{\sum j_i}}
\leq
\]
\[
\leq
\sum_{j_1=0}^{l_1-1}  \cdots  \sum_{j_m=0}^{l_m-1} 
\left(
\prod_{i=1}^m
{l_i-1\choose j_i}(l_i-2)^{{l_i-1-j_i}}
\left({1\over m} \right)^{j_i}
\right)
n^{{\sum j_i}}
=
\]
\[
  =\prod_{i=1}^{m}
  \left(
    \sum_{j_i=0}^{l_i-1} 
       {l_i-1\choose j_i}(l_i-2)^{l_i-1-j_i}
       \left(
          {1\over m}
       \right)^{j_i}n^{j_i}
  \right) =
\prod_{i=1}^{m}\left(l_i-2+{n\over m}\right)^{l_i-1},
\]
as wanted.
\end{Proof}

\bigskip This lemma is crucial to prove the two theorems announced in Section 2.
Indeed, Theorem \ref{thm.expectedsize} is just a 
version of the following more precise statement:

\begin{lemma}
\label{lema.expectedsize}
Consider polytopes $P$ of dimension $l$ and $Q$ of dimension $n-1$. Let $m$
be such that 
$m\leq n$. Given a triangulation $T_0$ of $P\times\Delta^{m-1}$ of
weighted size $t_0$  and a triangulation $T_Q$ of $Q$, the expected size of
the multi-staircase
triangulation $T_{P\times Q}$ of $P\times Q$ is bounded above by
\[
|T_Q|t_0\left({n\over{m}}+l\right)^l.
\]
\end{lemma}

\begin{proof}
Lemma \ref{lema.expectation} implies that
\[ E\left(\prod_{i=1}^{m}{k_i-1+l_i-1 \choose l_i-1}\right)=
E\left(\prod_{i=1}^{m}{{(k_i-1+l_i-1)^{\underline{l_i-1}}} 
\over{(l_i-1)!}}\right)\leq
{1\over{\prod_{i=1}^{m}(l_i-1)!}}\left(l+{n\over m}\right)^l.
\]
Hence:
\[ E\left(\sum_{\sigma\in T_Q}\sum_{\tau\in
T_0}\prod_{i=1}^{m}{{k_i-1+l_i-1}\choose{l_i-1}}\right)\leq
\sum_{\sigma\in T_Q}\sum_{\tau\in
T_0}{1\over{\prod_{i=1}^{m}(l_i-1)!}}\left({{n}\over{m}}+l\right)^l=
\]
\[ =|T_Q|\left(\sum_{\tau\in
T_0}{1\over{\prod_{i=1}^{m}(l_i-1)!}}\right)\left({n\over{m}}+l\right)^l=
|T_Q|t_0\left({n\over{m}}+l\right)^l.\]
\end{proof}

\medskip 

\begin{proof}[Proof of Theorem \ref{thm.triangulaciones}] 
Let $t_0= {\rho_0^l m^l}$ be the weighted size of the 
triangulation in the statement.
For any $n\geq {lm\rho_0\over \epsilon}$:
\[
\rho_0+\epsilon \geq \rho_0 {m\over n}\left({n\over m} + l\right),
\]
\[
n^l\left(\rho_0+\epsilon\right)^l 
   \geq {\rho_0^l m^l}\left({n\over m} + l\right)^l
   = t_0\left({n\over m} + l\right)^l,
\]
\[
{(n+l-1)!\over (n-1)!}\left(\rho_0+\epsilon\right)^l 
   \geq t_0\left({n\over m} + l\right)^l.
\]

Theorem \ref{thm.expectedsize}, with $Q=I^{n-1}$ and $P=I^l$, tells us that
\[
{\phi_{n+l-1}}\le \phi_{n-1}t_0\left({n\over m} + l\right)^l,
\]
or, in other words,
\[
{(n+l-1)!{\rho_{n+l-1}}^{n+l-1} \over (n-1)!{\rho_{n-1}}^{n-1}}\le 
t_0\left({n\over m} + l\right)^l.
\]

So, we have proved the first part of the theorem:
\[
\forall \epsilon >0, \quad
\forall n\geq {lm\rho_0\over \epsilon},
   \quad {\rho_{n+l-1}}^{n+l-1}\leq{\rho_{n-1}}^{n-1}{(\rho_{0}+\epsilon)}^{l}.
\]

The second part of the statement follows from the first one
with arguments similar to those of Corollary \ref{cor.Haiman}.
Recursively, we have that:
\[
\forall \epsilon >0,\quad \forall i\in\N,
   \quad \forall n\geq {lm\rho_0\over \epsilon},
   \quad 
    {\rho_{n+il}}\leq{\rho_{n}}^{n\over n+il}
    {(\rho_{0}+\epsilon)}^{il\over n+il}.
\]
This implies that for the given $l$ and any fixed $n$, taking $\epsilon={lm\rho_0\over n}$ we get
\[
\lim_{i\to\infty} \rho_{n+il} \le
\rho_0+{lm\rho_0\over n}.
\]
That is, the sequence of indices congruent to $n$ modulo $l$ has limit bounded by the right hand side. Since we can make $n$ as big as we want and the rest of the right-hand side are constants, 
 the $l$ subsequences of indices
modulo $l$ have limit bounded by $\rho_0$,
hence the limit of the whole sequence has this bound.
\end{proof}

\section{Interpretation of our method via the Cayley Trick}
\label{sec.cayley}

The Cayley Trick allows to study triangulations of a product $P\times
\Delta^{n-1}$ as mixed subdivisions of the Minkowski sum $P+\cdots+P$
($n$ summands). We overview here this method, but the reader should
look at \cite{Hu-Ra-Sa} for more details.

Let $Q_1,\ldots,Q_m\,\,\subset\R^d$ be convex polytopes  of vertex sets $\A_i$.
Consider their Minkowski sum, defined as
\begin{displaymath}
\sum_{i=1}^{m}Q_i = \left\{ \,x_1+\cdots+x_m : x_i\in Q_i
\right\}\mbox{.}
\end{displaymath}

We understand $\sum_{i=1}^mQ_i$ as a {\em marked polytope}, whose associated
point configuration is $\sum_{i=1}^m
\A_i:=\{q_1+\cdots +q_m : q_i\in
\A_i\}$. Here a {\it marked polytope} is a pair
$(P,\A)$ where $P$ is a polytope and $\A$ is a finite set of points of
$P$ including all the vertices. {\it Subdivisions} of a marked polytope are
defined in \cite[Chapter 7]{GKZbook} (sometimes they are called subdivisions of
$\A$). Roughly speaking, they are the polyhedral subdivisions of $P$ which use
only elements of $\A$ as vertices  (but perhaps not all of them). They form a
poset under the refinement  relation. The minimal elements are the
triangulations
of
$\A$.

A subset $B$ of $\sum_{i=1}^m \A_i$ is called {\it mixed} if
$B=B_1+\cdots +B_m$ for some non-empty subsets $B_i\subset \A_i$,
$i=1,\ldots,m$.  A {\it mixed subdivision} of $\sum_{i=1}^mQ_i$ is a
subdivision of it whose cells are all mixed. Mixed subdivisions form a subposet
of the poset
of all subdivisions, whose minimal elements are called {\em fine mixed}, in
which every mixed cell is \emph{fine}, i.e., does not properly contain any
other mixed cell.

We call {\it Cayley embedding} of $\{Q_1,\ldots,Q_m\}$ the marked polytope
$$(\C(Q_1,\ldots,Q_m),\C(\A_1,\ldots,\A_m)) \mbox{ in }
\R^d\times\R^{m-1}$$  defined as follows: let $e_1,\ldots,e_m$ denote an affine
basis in $\R^{m-1}$ and $\mu_i:\R^d\to\R^d\times\R^{m-1}$ be the
inclusion given
by $\mu_i(x)=(x,e_i)$. Then we define
\begin{displaymath}
   \C(\A_1,\ldots,\A_m):=\cup_{i=1}^m \mu_i(\A_i), \qquad
   \C(Q_1,\ldots,Q_m):=\hbox{conv}(\C(\A_1,\ldots,\A_m))
\end{displaymath}

Each $Q_i$ is naturally embedded as a face in
$\C(Q_1,\ldots,Q_m)$. Moreover, the vertex set of $\C(Q_1,\ldots,Q_m)$
is the disjoint union of the vertices of all the $Q_i$'s. This induces
the following bijection between cells in $\C(Q_1,\ldots,Q_m)$ and
mixed cells in $Q_1+\cdots+Q_m$: To each mixed cell $B_1+\cdots +B_m$
we associate the disjoint union $B_1\cup\cdots\cup B_m$. To a cell $B$
in  $\C(Q_1,\ldots,Q_m)$, we associate the Minkowski sum $B_1+\cdots
+B_m$, where $B_i=B\cap Q_i$.

\begin{thm}[The Cayley Trick,\cite{Hu-Ra-Sa}]
   \label{thm.cayley} Let $Q_1,\ldots,Q_m\,\,\subset\R^d$ be convex
   polytopes.  The bijection just exhibited induces an
   isomorphism between the poset of all subdivisions of \
   $\C(Q_1,\ldots,Q_m)$ and the poset of mixed subdivisions of
   $\sum_{i=1}^mQ_i$.  In this isomorphism triangulations correspond to
   fine mixed subdivisions.
\end{thm}

\begin{remark}
\label{rem.cayley.P} 
With the previous definitions, for any polytope $P$:
\[
 P\times\Delta^{m-1}=\C(P,\stackldots{m},P).
\]
In particular, in our context the Cayley Trick provides the following
bijections between triangulations and mixed subdivisions:

\[
\begin{array}{ccc}
\framebox{\mbox{\protect\small Triangulation of
$P\times\Delta^{m-1}$}} & & \framebox{\mbox{\protect\small Triangulation of
$P\times\Delta^{n-1}$}}\\
  & & \\
\updownarrow \mbox{\protect\small Cayley Trick} & & \updownarrow
\mbox{\protect\small Cayley Trick}\\
  & & \\
\framebox{\mbox{\protect\small Fine mixed subdivision of
$P+ \stackcdots{m}+P$}} &
\longrightarrow & \framebox{\mbox{\protect\small Fine mixed subdivision of
$P+\stackcdots{n}+P$}}
\end{array}
\]
\end{remark}

Our interest in the Cayley Trick is two-fold. On the one hand, it
provides a way to visualize our candidate triangulations of
$I^l\times\Delta^{m-1}$ as objects in dimension $l$, instead of
$l+m-1$. We will use this in Section \ref{sec.rholm}. 

On the other
hand, the construction of the previous sections has a simple geometric
interpretation in terms of the Cayley Trick. More precisely, the
polyhedral subdivision of Corollary \ref{coro.subdivision} can be obtained as
follows: Let $S$ be a polyhedral subdivision of $P\times
\Delta^{m-1}$, and $S_M$ the corresponding mixed subdivision of 
$P+ \stackcdots{m}+P$. Each cell in $S$ decomposes uniquely as
\[
\cup_{i=1}^m B_i\times\{v_i\},
\]
where the $B_i$ are subsets of vertices of $P$. The corresponding cell
in $S_M$ is just $B_1+\cdots+B_m$. To construct our polyhedral subdivision
of $P\times \Delta^{n-1}$ we just need to scale each summand $B_i$ by
the integer $k_i$ which tells us how many vertices of $\Delta^{n-1}$
correspond to the vertex $v_i$ of $\Delta^{m-1}$.

That is to say, from $S_M$ we construct a mixed subdivision of 
$P+ \stackcdots{n}+P$ by the formula
\[
\widetilde{S_M}=
 \{B_1+ \stackcdots{k_1}+B_1+\cdots \cdots +B_m+ \stackcdots{k_m}+B_m
 : B_1+\cdots+B_m \in S_M\}.
\]
The polyhedral subdivision $\tilde S$ of $P\times\Delta^{n-1}$ stated in
Corollary \ref{coro.subdivision} is the one corresponding via the
Cayley Trick to the mixed subdivision $\widetilde{S_M}$ of $P+
\stackcdots{n}+P$.

Also the type and weight of a simplex in $P\times\Delta^{m-1}$
have a simple interpretation via the Cayley Trick. With the
notation of Definition \ref{def.main}, let $\tau=\tau_1\cup\cdots\cup\tau_m$
be a simplex. The corresponding cell $\tau_M$ in $P+\stackcdots{m}+P$ is the
Minkowski sum of the simplices $\tau_1,\dots,\tau_m$, which lie in
complementary affine subspaces. Hence, $\tau_M$ is combinatorially a
product of $m$ simplices, of dimensions $t_1,\dots,t_m$ where 
$(t_1,\dots,t_m)$ is the type of $\tau$. Then the weight of $\tau$
represents the volume of $\tau_M$,
normalized with respect to the unit parallelepiped in the lattice spanned by
the vertices of $\tau_M$. With this we can prove:

\begin{prop}
\label{prop.weightedsize}
Let $P$ be a lattice polytope of dimension $l$. Let $V$ be its volume, 
normalized to the unit
parallelepiped in the lattice. Then, the weighted size of a
triangulation of $P\times\Delta^{m-1}$ is at most $m^l V$, with equality
if and only if the triangulation is unimodular (with respect to the
lattice). 

In particular, the weighted efficiency of a triangulation of
$I^l\times \Delta^{m-1}$ is at most one, with equality for unimodular
triangulations.
\end{prop}

\begin{Proof}
If $P$ is a lattice polytope (e.g., a cube), then
the lattice spanned by $\tau_M$ is a sublattice of the one spanned
by the point configuration $P+\stackcdots{m}+P$, and coincides with it
if and only if $\tau$ is unimodular. 
Then, for unimodular triangulations
the weighted size is then just the volume of $P+\stackcdots{m}+P$,
normalized to the unit parallelepiped, which equals $m^l V$. For
non-unimodular triangulations the weighted size is smaller than that.
\end{Proof}

\section{More on $\rho_{l,m}$}
\label{sec.rholm}

In this section we obtain the value of $\rho_{2,m}$ for any $m$ and
we will also show triangulations of $I^3\times\Delta^1$ and 
$I^3\times\Delta^2$ providing the values of
$\rho_{3,2}$ and $\rho_{3,3}$ stated in Section
\ref{rhos.cplex}.

\subsection{Smallest weighted efficiency of triangulations of
$I^2\times\Delta^{m-1}$}
\label{subsec.l=2}

\begin{figure}[htb]
\begin{center}
         \epsfxsize=2.5 in\leavevmode \epsfbox{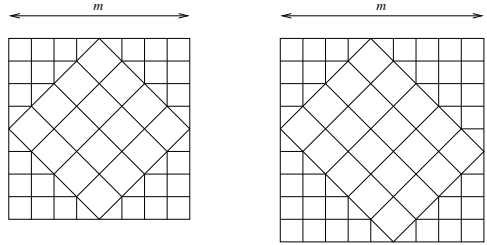}
\end{center}
\caption{Fine mixed subdivisions corresponding to triangulations of
$I^2\times\Delta^{m-1}$ with smallest
  weighted efficiency for $m$ even (left) or odd (right).}
\label{fig.m=2}
\end{figure}

Here we prove:

\begin{thm}
\label{thm.m=2}
The smallest weighted efficiency $\rho_{2,m}$ of
triangulations of $I^2\times\Delta^{m-1}$ is
\[
\rho_{2,m}=\sqrt{{\lceil {3m^2}/4 \rceil}\over{m^2}}.
\] 
That is, $\sqrt{3/4}$ for even $m$ and $\sqrt{3/4}+\Theta(m^{-2})$ for odd $m$.
A fine mixed subdivision of \mbox{$I^2+\stackcdots{m}+I^2$} corresponding to a
triangulation of $I^2\times\Delta^{m-1}$ with that weighted
efficiency is given in Figure \ref{fig.m=2}.
\end{thm}

Let $B_M=B_1+\stackcdots{m}+B_m$ be a cell in a fine mixed subdivision of 
$I^2+\stackcdots{m}+I^2$ (a square of size $m$). Each $B_i$ must be a
simplex coming from the $i$-th copy of $I^2$, and in order to have a
{\em fine} mixed subdivision, the different $B_i$'s must lie in
complementary affine subspaces. Then, there are the following
possibilities:

\begin{itemize}

\item $B_M$ is a triangle, obtained as the sum of a triangle in one of
the $I^l$'s and a single point in the others. The weight of $B_M$ is $1/2$.

\item $B_M$ is a quadrangle, obtained as the sum of two (non-parallel)
segments from
two of the $B_i$'s and a point in the rest of them. Three types of
quadrangles can appear, depending on whether none, one or both of the
two segments involved is a diagonal of $I^2$: 
a square parallel to the axes, a rhomboid (both with area 1) and
a diagonal square (of area 2).
The weight of $B_M$ is $1$.

\end{itemize}

In particular, the weighted size of the mixed
subdivision equals $T/2 + S_1 + S_2 + R$, where $T$, $S_1$, $S_2$ and $R$
denote, respectively, the numbers of triangles, squares of area 1, squares of
area 2 and rhombi in the subdivision. Since the total area of
$I^2+\stackcdots{m}+I^2$ is $m^2 = T/2 + S_1 + 2S_2 +R$, we conclude
that:

\begin{prop}
\label{prop.squares}
The weighted size of a mixed subdivision of $I^2+\stackcdots{m}+I^2$
equals $m^2 - S_2$, where $S_2$ is the number of squares of size 2
in the subdivision.
\end{prop}

With this we can already conclude that the fine
mixed subdivisions shown in Figure \ref{fig.m=2} (one
for $m$ even and one for $m$ odd) have weighted size equal to
$\lceil 3m^2/4 \rceil$, since they have exactly $\lfloor
m^2/4\rfloor$ squares of area 2.

Our task is now to prove that no mixed subdivision can have more than 
$\lfloor m^2/4\rfloor$ squares of area 2. For this we use:

\begin{lemma}
\label{lemma.mixtopoly}
Taking the $i$-th summands of all the mixed cells
$B_1+\stackcdots{m}+B_m$ in a mixed subdivision of
$P_1+\stackcdots{m}+P_m$ produces a polyhedral subdivision of
$P_i$. If the mixed subdivision was fine, the polyhedral subdivision
is a triangulation.
\end{lemma}

\begin{Proof} 
By the Cayley Trick, the mixed subdivision of $P_1+\stackcdots{m}+P_m$
induces a polyhedral subdivision of  $\C(P_1,\dots,P_m)$, and a
triangulation if the mixed subdivision is fine. Since 
the polytope $P_i$ appears as a face in $\C(P_1,\dots,P_m)$, 
any subdivision (resp. triangulation) of
$\C(P_1,\dots,P_m)$ induces a subdivision (resp. triangulation) of
$P_i$. That this subdivision is the one obtained taking the $i$-th
summands of all the mixed cells follows from Theorem \ref{thm.cayley}.
\end{Proof}

\begin{prop}
\label{prop.l=2}
Let $S_M$ be a fine mixed subdivision of $I^2 +\stackcdots{m} +
I^2$. Let $a$ and $b$ be the number of summands $I^2$'s which are
triangulated in one and the other possible triangulations of $I^2$
(i.e., using one or the other diagonal). Then, the number of squares
of area 2 in $S_M$ is at most $ab$ and the weighted size of $S_M$ is
at least $m^2-ab$.
\end{prop}

\begin{Proof}
Each square of area 2 is the Minkowski sum of two opposite diagonals
of two copies of $I^2$ (say, the $i$th and $j$th copies), 
and a point in the other $m-2$ of the copies
of $I^2$. Clearly, the two copies which contribute diagonals have to
be triangulated in opposite ways. The only thing which remains to be
shown is that the same pair of copies of $I^2$ cannot contribute two
different squares of size 2. For this, observe that ``contracting'' in
every mixed cell of a mixed subdivision all the summands other than
the $i$th and $j$th should give a mixed subdivision of $P_i+P_j$. And
in a mixed subdivision of two squares there is no room to put two
different diagonal squares of area 2.
\end{Proof}

In the statement of Proposition \ref{prop.l=2}, we have that $a+b =m$. In
particular, the maximum possible value of $ab$ is $\lfloor
m^2/4\rfloor$. This finishes the proof of Theorem \ref{thm.m=2}.

\subsection{Smallest weighted efficiency triangulations of
$I^{3}\times\Delta^{m-1}, m=2,3$}
\label{subsec.visual}

Here we will try to visualize triangulations of
$I^{3}\times\Delta^{1}$ and $I^{3}\times\Delta^{2}$ which minimize 
the weighted efficiency, which we computed using the integer
programming approach sketched in Section \ref{sec.overview}.

Of course, we use the Cayley trick to decrease the dimension, so we
show the corresponding fine mixed subdivision of a Minkowski sum instead of
the triangulation itself.

\begin{itemize}
\item \framebox{$I^3\times\Delta^{1}$} \qquad \qquad
$I^3\times\Delta^{1}=\C(I^3,I^3)\longleftrightarrow I^3+I^3$
\smallskip

We have to give a subdivision of a 3-cube of size 2. For this we first
cut the eight corners of the cube, producing a \emph{cubeoctahedron}, a
semi-regular 3-polytope with 6 square and 8 triangular facets. The
edges of the cubeoctahedron can be distributed in four ``equatorial
hexagons'' each of which cuts the polytope into 2 halves. 
Figure \ref{fig.I3D1} depicts these two halves for one of the
equatorial hexagons. The labels in the vertices are heights,
interpreted as follows: Our point configuration is $\{0,1,2\}^3$ and
the height of point $(i,j,k)$ is just $i+j+k$.

It turns out that performing three of these four possible
halvings, the cubeoctahedron is
decomposed into six triangular prisms and two tetrahedra. In Figure
\ref{fig.I3D1} each half is actually decomposed into three prisms and
one tetrahedron. These,
together with the eight tetrahedra we have 
cut from corners, form a fine mixed
subdivision of $\{0,1,2\}^3$ with 10 tetrahedra and 6 triangular prisms. Its 
weighted size is then
\[
10 \cdot {1\over 6} + 6 \cdot{1\over 2} = {14\over 3}
\]
and its weighted efficiency $\sqrt[3]{{14/3}\over 8}$.

\smallskip

\begin{figure}[htb]
\begin{center}
         \epsfxsize=3 in\leavevmode \epsfbox{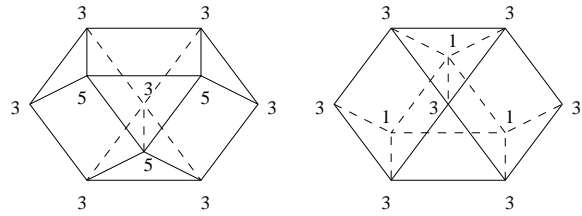}
\end{center}
\caption{Fine mixed subdivision of $I^3+I^3$ corresponding to a
triangulation of
$I^3\times\Delta^{1}$ with smallest weighted efficiency.}
\label{fig.I3D1}
\end{figure}

\item \framebox{$I^3\times\Delta^{2}$} \qquad \quad
$I^3\times\Delta^{2}=\C(I^3,I^3,I^3)\longleftrightarrow I^3+I^3+I^3$
\smallskip

The fine mixed subdivision of $I^3+I^3+I^3$ is given in Figure \ref{fig.I3D2}. It 
consists of 20 triangular prisms, 16 tetrahedra and 2 parallelepipeds. 

Thus, the weighted size of the triangulation is $20{1\over 2}+16{1\over
6}+2{1\over 1}={44\over 3}$, and then the smallest weighted efficiency is
\[
\rho_{3,3}=\sqrt[3]{44/3\over 3^3}.
\]

\smallskip

\begin{figure}[htb]
\begin{center}
         \epsfxsize=5 in\leavevmode \epsfbox{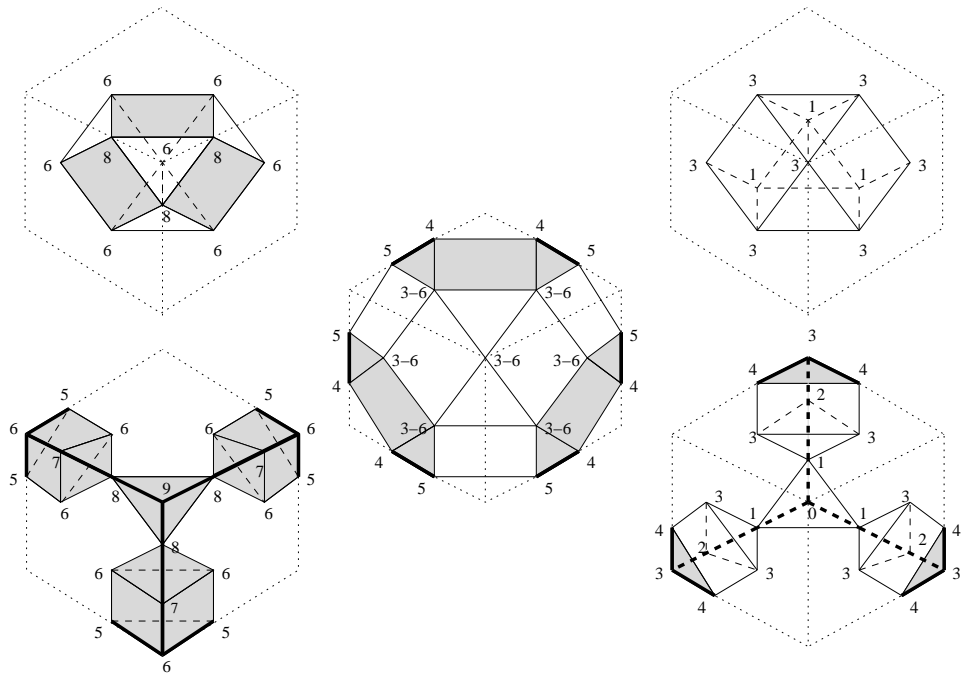}
\end{center}
\caption{Fine mixed subdivision of $I^3+I^3+I^3$ corresponding to a
triangulation
of $I^3\times\Delta^{2}$ with smallest weighted efficiency.}
\label{fig.I3D2}
\end{figure}

\end{itemize}

Let us explain how to interpret  Figure \ref{fig.I3D2}. Again, each point $(i,j,k)$ in
$\{0,1,2,3\}^3$ has been given a height $i+j+k$, which is written next
to it. The subdivision is displayed in five parts. 
The left-top portion in the figure
is a half cubeoctahedron
exactly as the one in Figure \ref{fig.I3D2}. The reader has to assume it divided
into a tetrahedron and three triangular prisms, as before.
The left-bottom portion consists of a corner
tetrahedron and three triangular prisms located at corners of the cube, each
joined to the half cubeoctahedron by a tetrahedron. So far we have six prisms and
five tetrahedra, and the same is got from the right-top and right-bottom portions of the Figure.

In between the
two half cubeoctahedra, however, we have now a hexagonal prism decomposed into
two triangular prisms and two quadrilateral prisms, all of height
$\sqrt{3}$. The hexagonal prism is surrounded by a belt formed by six
triangular prisms and six tetrahedra.

The thick edges in Figure \ref{fig.I3D2} represent edges of the big 3-cube of
side 3, whose three visible facets are drawn by dotted lines. 
We have also shaded the facets of the subdivision which are contained in
those facets of the big 3-cube.

%%%%%%%%%%%%%%%%%%%%%%%%%%%%%

\end{document}